\documentclass[dvips, 12pt,a4paper]{article}

\usepackage[greek,american]{babel}
\usepackage{amssymb}
\usepackage{amsfonts}
\usepackage{amsmath}
\usepackage{subfigure}
\usepackage[dvips]{graphicx}

\begin{document}

\begin{center}
{{{\Large\sf\bf Bayesian Nonparametric Density Estimation under Length Bias }}}\\

\vspace{0.5cm}
{\large\sf Spyridon J. Hatjispyros
$^*$, Theodoros Nicoleris$^{**}$
and Stephen G. Walker$^{***}$ }\\

\vspace{0.2cm}
\end{center}
\centerline{\sf $^{*}$ Department of Mathematics, University of the Aegean,}
\centerline {\sf Karlovassi, Samos, GR-832 00, Greece.} 
\centerline{\sf $^{**}$ Department of Economics, National and Kapodistrian University of Athens,}
\centerline{\sf Athens, GR-105 59, Greece. }  
\centerline{\sf $^{***}$Department of Mathematics, University of Texas at Austin,}
\centerline{\sf Austin, Texas 7812, USA. }

\begin{abstract}

A density estimation method  in a  Bayesian nonparametric framework is presented
when recorded data are not  coming directly from the distribution of interest, but from a  length biased version.
From a Bayesian perspective,  efforts to computationally evaluate   posterior quantities conditionally on length biased data   
were hindered  by the inability   to circumvent the problem of a normalizing constant. In this paper we present a novel 
Bayesian nonparametric approach to the length bias sampling problem which circumvents the issue of the normalizing constant. 
Numerical illustrations as well as a  real data example are presented and the estimator is compared against its frequentist 
counterpart, the kernel density estimator for indirect data of Jones (1991).

\vspace{0.1in} \noindent {\sl Keywords:} Bayesian nonparametric inference; Length biased sampling; Metropolis algorithm.

\end{abstract}

\section{Introduction}

Let $f(x;\theta)$, with $\theta\in\Theta$ being an unknown parameter, be a family of density functions.
Sampling under selection bias involves observations being drawn not from $f(x;\theta)$ directly, 
but rather from  a distribution which is a biased  version   of $f(x;\theta)$, given by the density function
$$g(y;\theta) =\frac{ w(y)\,f(y;\theta)}{\int_0^\infty w(x)\,f(x;\theta)\,d x }.
$$
where the $w(x)>0$ is the weight function. We observe a sample  $(Y_1,\ldots, Y_n )$, independently taken from $g(\cdot)$. 
In particular, when the weight function is linear; i.e.  $w(y)=y$, the samples are known as length biased.

There are many  situations where weighted data arise; for example, in survival analysis (Asgharian et al., 2002); 
quality control problems for estimating fiber length distributions (Cox, 1969); models with clustered  or over--dispersed 
data (Efron, 1986); visibility bias  in aerial data;  sampling from queues or telephone networks. For further examples 
of length biased sampling see, for example, Patil and Rao (1978) and Patil (2002). 

In the nonparametric setting  $f(x;\theta)$ is replaced by the more general $f(x)$, so  the likelihood function for $n$ data points becomes,
$$l(f)=\prod_{i=1}^n \frac{y_i\,f(y_i)}{\int_0^\infty x\,f(x)\,d x}.$$
A classical nonparametric maximum likelihood estimator (NPMLE) for $F$ (the disribution function corresponding to $f$)  exists for this problem and is discrete, with atoms located at the observed data points. In particular, Vardi (1982) finds an explicit form for the estimator in the presence of two independent samples, one  from $f$ and the other  from  the   length biased density $g$.

Our work focuses on length biased sampling  and from the  Bayesian nonparametric setting we work in, the aim is to obtain a density estimator for $f$. 
There has been no work done on this problem in the Bayesian nonparametric framework due to the issue of the intractable likelihood function, particularly when $f$ is modeled nonparametrically using, for example, the mixture of Dirichlet process (MDP) model; see Lo (1984) and Escobar and West (1995). While some ideas do exist on how to deal  with intractable normalizing constants; see Murray et al. (2006); Tokdar, (2007); Adams et al. (2009); and Walker, (2011), these ideas fail here for two reasons: the infinite dimensional model and the unbounded $w(y)=y$ when the space of observations is the positive reals.     

We by-pass the intractable normalizing constant by modeling $g$ nonparametrically. We argue that modeling $f$ or $g$ nonparametrically is providing the same flexibility to either; i.e. modeling $f(y)$ nonparametrically and defining $g(y)\propto y\,f(y)$ is essentially equivalent to modeling $g(y)$ nonparametrically and defining $f(y)\propto y^{-1}g(y)$. We adopt the latter style, obtain samples from the predictive density of $g$ and then ``convert" these samples from $g$ into samples from $f$, which forms the basis of the density estimator of $f$.

The layout of the paper is as follows:  In Section 2 we provide supporting theory for the model idea which avoids the need to deal with the  intractable likelihood function. Section 3 describes the model and the MCMC algorithm for estimating it and Section 4 describes some numerical illustrations. In Section 5 are  the concluding remarks and in Section 6 asymptotic results are provided. 

\section{Supporting theory and methodology}

Our aim is to avoid computing the intractable normalizing constant. The strategy for that would be  to model the density $g(y)$ directly and then make inference about $f(x)$ by exploiting the fact that
$$
g(y)\propto y\,f(y).
$$

 In the parametric case if a  family $f(x;\theta)$ is known  then so is $g(y;\theta)$, except its normalizing constant may not be tractable. There is a reluctance to avoid the problem of the normalizing constant in the parametric case by modeling the data directly with a tractable $g(y;\theta)$ since the incorrect model would be employed. However, in the nonparametric setting it is not regarded as relevant whether one models $f(\cdot)$ or $g(\cdot)$ directly. A clear motivation to model  $g(\cdot)$ directly is that this is where the data are coming from. 

For  a general weight function $w$,  an essential condition  to model $F$  through $G$ ($F$ and $G$ denote the corresponding distribution functions of $f$ and $g$, respectively)  is the finiteness of   $\int_0^\infty w(x)^{-1}\,g(x)\,d x$. This, through invertibility,  enables us to reconstruct $F$ from $G$ and occurs when $F$ is  absolutely continuous with respect to $G$, with the Radon-Nikodym derivative being  proportional to $w(x)^{-1} $.

For absolute continuity to hold we need that $w(x)>0$ in the support of  $F$ ie $F(x:w(x)=0)=0$. In the length biased case examined here $w(x)=x$ and the densities have support  on the positive real line, so this condition is automatically satisfied. A  case, for instance, when this does not hold  and invertibility fails is in a truncated model where $w(x)=1(x\in B)$, $B$ is a Borel set and $F$ is a distribution  which  could be  positive outside of $B$.

A  Bayesian model is thus  constructed by assigning an appropriate nonparametric prior  distribution   to $g$, provided that
$$\int_0^\infty y^{-1}\,g(y)\,d y<\infty.$$
This in turn specifies a prior for $f$.

The question that now arises  is how the posterior structures obtained after modelling $g$ directly can be converted to posterior structures from $f$.
The first step in this process would be to devise a method  to convert a  biased sample from a   density $g$  to one  from its  debiased version  $f$. This  algorithm is then incorporated to  our model building process so that posterior inference becomes  possible.
 
Specifically, assume  that a sample $y_1,\ldots,y_N$, comes from a biased density  $g$. This  can be  converted into a sample from  $f(x)\propto x^{-1}g(x)$  using a Metropolis--Hastings algorithm. If we denote  the current sample from $f(x)$ as $x_j$, then
$$x_{j+1}=y_{j+1}\quad\mbox{with probability}\quad\min\left\{1,\frac{x_j}{y_{j+1}}\right\},$$
otherwise $x_{j+1}=x_j$.  Here, we have the transition density for this process as
$$p(x_{j+1}|x_j)=\min\left\{1,\frac{x_j}{x_{j+1}}\right\}g(x_{j+1})+\left\{1-r(x_j)\right\}{\bf 1}(x_{j+1}=x_j),$$
where
$$r(x)=\int \min\left\{1,\frac{x}{x^*}\right\}g(x^*)\,d x^*.$$
This transition density satisfies detailed balance with respect to  $f(x)$ since
$$p(x_{j+1}|x_j)\,x_j^{-1}g(x_j)=p(x_j|x_{j+1})\,x_{j+1}^{-1}g(x_{j+1})$$
and thus the transition density has stationary density given by $f(x)$.

This algorithm was  first tested  on    a toy   example, i.e.   $g(y)$ is  Ga$(y|2,1)$ so that $f(x)$ is Ga$(x|1,1)$. A sample of $N=10,000$ of the  $(y_j)$ was taken independently from the $g(\cdot)$ and the Metropolis algorithm run to generate the $(x_j)$, starting with $x_0=1$. Sample values for the sequence of $(x_j)$ yield
$$N^{-1}\sum_{j=1}^N x_j=0.981\quad\mbox{and}\quad N^{-1}\sum_{j=1}^N x_j^2=1.994,$$
which are compatible outcomes with the $(x_j)$ sample coming from $f(x)$. A similar example will be elaborated on in the numerical illustration section.

Applying this idea to   our model would amount to turning a sample from the biased posterior predictive density  to an unbiased one using a MH  step. An outline  of the inferential methodology is now described.  
\begin{enumerate}
\item Once data  $(y_1,\ldots,y_n)$  from a biased distribution $g$ become avalaible a model for $g$ is assumed and a nonparametric prior is assigned. 
\item Using MCMC methods, after a sensible burn-in period, at each iteration, posterior  values of the random measure $\Pi(d g|y_1,\ldots,y_n)$ and the relevant parameters are obtained. Subsequently, conditionally on those values, a sequence  $\{y_{n+1}^{(l)}\},\;l=1,2,\ldots$,  from the posterior predictive density $g(y|y_1,\ldots,y_n)$ is generated.
\item The $\{y_{n+1}^{(l)}\}$ will then  form a sequence of  proposal values of a Metropolis-Hastings chain  with stationary density the debiased version of  the posterior predictive, i.e. $\propto y^{-1}g(y|y_1,\ldots,y_n)$. 
 Specifically, at the $l$-th iteration of the algorithm  applying  a rejection  criterion a value   $x_{n+1}^{(l)}$ is  generated such that $x_{n+1}^{(l)}=y_{n+1}^{(l)}$  with probability $\min\left\{1,x_{n+1}^{(l-1)}/y_{n+1}^{(l)}\right\}$, otherwise $x_{n+1}^{(l)}=x_{n+1}^{(l-1)}$. 
\item These $\{x_{n+1}^{(l)}\}$ values form a sample  from the posterior predictive of $f$.
\end{enumerate}

\section{ The model and inference}

We want the model for $g(y)$ to have large support and the standard Bayesian nonparametric idea for achieving this is based on infinite mixture models (Lo, 1984) of the type
$$g_{\mathbb P}(y)=\int k(y;\theta)\,{\mathbb P}(d\theta)$$
where ${\mathbb P}$ is a discrete probability measure and $k(y;\theta)$ is a density on $(0,\infty)$ for all $\theta$. Since we require $g(y)$ to be such that
$$\int_0^\infty y^{-1}\,g_{\mathbb P }(y)\,d y<\infty $$ or, equivalently, for a kernel $k(y;\theta)$
$$\int_0^\infty \,y^{-1}k(y;\theta)dy<\infty$$
we find it most appropriate to take the kernel to be a log--normal distribution. So, assuming a constant precision parameter $\lambda$ for each component, we have
\begin{equation}
\label{randomdensity1}
g_{\lambda,{\mathbb P}\,}(y)=
\int_{\mathbb R} \mbox{LN}(y\vert\,\mu,\lambda^{-1})\,{\mathbb P}(d\mu)
\end{equation}
where ${\mathbb P}$ is a discrete random probability measure  defined in
${\mathbb R}$ and ${\mathbb P}\sim {\rm DP}(c,P_0)$, where ${\rm DP}(c,P_0)$
denotes the Dirichlet process (Ferguson, 1973) with precision parameter
$c>0$ and base measure $P_0$. Interpreting the parameters, we have that
${\rm E}( {\mathbb P(A)})=P_0(A)$, and
$$
{\rm Var}({\mathbb P}(A))={P_0(A)(1-P_0(A))\over c+1}
$$
for appropriate sets $A$.

This  Dirichlet process mixture model implies the 
hierarchical model for $y=(y_1,\ldots,y_n)$:
For $1\le i\le n$
\begin{eqnarray}
\nonumber
& & y_i\vert\, \mu_i, \lambda \,\stackrel{\rm ind}\sim\, \mbox{LN}(\mu_i,\lambda^{-1})\\
\nonumber
& & \mu_i\vert\, {\mathbb P}\,\stackrel{\rm i.i.d.}\sim\, {\mathbb P}\\
\nonumber
& & {\mathbb P}\vert\,c,P_0\,\sim\, {\rm DP}(c,P_0)
\end{eqnarray}
To complete the model we choose  $\lambda\sim$ Ga$(a,b)$ and for the base measure, $P_0$ is ${\rm N}(0,s^{-1})$.

A useful representation of the Dirichlet process, introduced  by Sethuraman and Tiwari (1982)
and Sethuraman (1994),
is the stick--breaking constructive representation given by
$$
{\mathbb P}\,=\,\sum_{j=1}^\infty w_j\,\delta_{\mu_j}
$$
where the $(\mu_j)$ are i.i.d. from $P_0$, i.e.  ${\rm N}(0,s^{-1})$. The $(w_j)$ are constructed
via a stick--breaking process; so that $w_1=v_1$ and, for $j>1$,
\begin{equation}
\label{stick}
w_j=v_j\prod_{l<j}(1-v_l)
\end{equation}
where the $(v_j)$ are i.i.d. from the beta$(1,c)$ distribution, for some $c>0$,
and $\sum_{j=1}^\infty w_j=1$ almost surely. Let $w=(w_j)_{j=1}^\infty$ and
$\mu=(\mu_j)_{j=1}^\infty$; then we can then write
\begin{equation}
\label{randomdensity2b}
g(y_i\vert \mu,w,\lambda)=
\sum_{j=1}^\infty w_j\,\mbox{LN}\left(y_i\vert \mu_j,\lambda^{-1}\right)
\end{equation}
This is a standard  Bayesian nonparametric model. 
The MCMC algorithm is implemented  using latent variable techniques, despite the infinite dimensional model. 
The basis of this sampler is in Walker (2007) and Kalli et al. (2011).

For $1\le i\le n$  we introduce latent variables $u_i$ which make the sum finite.
The $u_i$ augmented  density then becomes,
\begin{equation}
\begin{split}
\label{randomdensity3} g(y_i, u_i\vert\, w, \mu, \lambda)\,= \,
\sum_{j=1}^\infty {\bf 1}(u_i<w_j)\mbox{LN}\left(y_i\vert\mu_j,\lambda^{-1}\right)\, = \\
=\sum_{j\in A_w(u_i)}\mbox{LN}\left(y_i\vert\mu_j,\lambda^{-1}\right)
\end{split}
\end{equation}
This   has a finite representation and  $A_w(u_i)$ denotes the
almost surely finite $u_i$ slice set $\{j:u_i<w_j\}$.

Now we introduce latent variables $\{d_1,\ldots,d_n\}$ which allocate  the  component
that $\{y_1,\ldots,y_n\}$ are sampled  from. Conditionally on the weights $w$
these are sampled  independently with  $P(d_i=j|w)=w_j$.
Hence,  we consider the $(u_i,d_i)$ augmented random density
$$
g(y_i,u_i,d_i\vert\,w, \mu, \lambda)\,=\,{\bf 1}(u_i<w_{d_i})
\mbox{LN}\left(y_i\vert\mu_{d_i},\lambda^{-1}\right)
$$
Therefore, the complete data likelihood based on a sample of size $n$  is seen to be
$$
l(y,u,d\vert\,w, \mu, \lambda)\,=\,
\prod_{i=1}^n {\bf 1}(u_i<w_{d_i})\mbox{LN}\left(y_i\vert\mu_{d_i},\lambda^{-1}\right)
$$
This will form the basis of our Gibbs sampler.
At each iteration we sample from the associated full conditional
densities of the following variables:
\begin{eqnarray}
\nonumber
  (v_j,\mu_j),\;j=1,\ldots,N;\;\lambda;\;(d_i,u_i),\;i=1\ldots n  & &
\end{eqnarray}
where $N$ is a random variable, such that
$\cup_{i=1}^n A_w(u_i)\subseteq \{1,\ldots,N\}$, and $N<\infty$ almost surely.

These distributions are, by now, quite standard so we proceed directly to the last two steps of the algorithm. 

The upshot is that after a sensible burn--in time period given  the current selection of parameters,  at each  iteration, we can  sample values $y_{n+1}$  from  the posterior predictive density  $g(y|y_1,\ldots,y_n)$ and  subsequently, using a Metropolis step, draw a $z$ value from its debiased version  $f(\cdot|y_1,\ldots,y_n)$.
\begin{itemize}
\item [{\bf 1.}]
Once  stationarity is reached then at each iteration we have points generated  by the  posterior measure  of the variables. These
points  are represented by
$$\{ (v_j^*, \mu_j^*),\;j=1,2\ldots\;;\lambda^*;\;(d_i^*,u_i^*),\;i=1\ldots
n;\;\}$$
Given  $\{ w_j^*,\mu^*_{j},\lambda^*\}$  a value $y_{n+1}\,\sim g(y|y_1,\ldots,y_n)$ is generated.                        
This is done by sampling  a $r$  uniformly in the unit interval and then take $k=1$ if $0<r\leq w^*_1$
or $k\geq 2$ if
$$
\sum_{i=1}^{k-1}w_j^* < r \leq \sum_{j=1}^{k}w_j^*
$$ 
The appropriate $\mu_{n+1}^*=\mu_{k}^*$ is then assigned, with probability $w^*_{k}$. Even though we have not sampled all the weights, if we ``run out" of weights, in essence the indices \{1,\ldots,N\},  we merely draw a $\mu_{n+1}^*$ from $ {\rm N}(0,s^{-1})$. Finally, the predictive value $y_{n+1}$  comes from  $\mbox{LN}(y \vert\mu_{n+1}^*,\lambda^*)$. 
\item [{\bf 2.}] The Metropolis step for the posterior predictive of $f$: Let $\tilde{x}$ be the state of the chain
from the previous Gibbs iteration. Accept the sample $y_{n+1}$, from the $g$-predictive, as coming from the $f$-predictive,
that is $z=y_{n+1}$, with probability $\min\{1, \tilde{x}/y_{n+1}\}$; otherwise the chain remains in its current state i.e. $z=\tilde{x}$.

\end{itemize} 

\section{Numerical illustrations}

We illustrate the model with two simulated data sets and a real data example. In each of the assumed models, 
for a given realisation $(y_1,\ldots,y_n)$, we report on the  results and compare  them  with the following density estimators:
\begin{itemize}
\item[{\bf (i)}] The classical kernel density estimate given by
\begin{equation}
\label{USUALKDE}
{\widetilde g}_h\,(y;\,(y_1,\ldots,y_n))\,\propto\,
n^{-1}\,\sum_{j\,=1}^n {\rm N}\left(y|\, y_j,h^2\right){\bf 1}(y>0).
\end{equation}
\item [{\bf (ii)}] The kernel density estimate for indirect data, see Jones (1991), is given by
\begin{equation}
\label{NPMLE}
{\widehat f}_{J,h}\,(y;\,(y_1,\ldots,y_n))\,\propto\,
n^{-1}{\widehat\mu}\,\sum_{j\,=1}^n y_j^{-1}\,{\rm N}\left(y|\, y_j,h^2\right){\bf 1}(y>0),
\end{equation}
where  ${\widehat\mu}$ is the harmonic mean of $(y_1,\ldots,y_n)$.
\end{itemize}
Here $h$ is the bandwidth and in all cases an estimate of it has been calculated 
as the average of the plug--in and solve--the--equation versions of it,
(Sheather and Jones $1991$).
The Gibbs sampler iterates $60,000$ times with a burn--in period of $10,000$.

\subsection{Simulated Data Examples}

Here we use  non informative prior specifications:  
\begin{equation}\label{noninformative}
\pi(\lambda)\propto 1/\lambda,\quad{\rm and}\quad\mu_j\sim{\rm N}(0,0.5^{-1}).
\end{equation}
The value of the concentration parameter has been set to $c=1$.

\vspace{0.1in}
\noindent
{\bf Example 1.} The length biased distribution is $g(y)={\rm Ga}(y|\,2,0.5)$ and we simulate
$y_{\rm g}=(y_1,\ldots,y_n)$ of size $n=50$. The following results are presented  Figure 1:

\begin{itemize}
\item 1(a): 
(i) a histogram  of the simulated length biased data set $y_{\rm g}$, ii) the true biased  
density Ga$(2,0.5)$ (the solid line) and iii)
the kernel density estimate  $\widetilde{g}_h(y;y_{\rm g})$ (the dashed line).
\item 1(b):  
(i) a histogram  of a sample from the posterior predictive density  $g(y_{n+1}|y_{\rm g})$,
(ii) the true biased  density Ga$(2,0.5)$ (the solid line) and iii) the kernel density estimate   $\widetilde{g}_h(y;y_{\rm g})$ (the dashed line).
\item 1(c): (i) a histogram  of the debiased data associated with the application of the Metropolis step, ii) the    true  debiased   density $\exp(0.5)$ (the solid line)   and iii) Jones' kernel density estimate ${\widehat f_{J,h}(y;y_{\rm g}})$ (the dashed line).
\end{itemize}
For both estimators $\widetilde{g}_h(y;y_{\rm g})$ and ${\widehat f}_{J,h}(y;y_{\rm g})$ the bandwidth parameter is set at 
$h=1.06$. The average number of clusters ${\cal C}_g$ is  $4.27$. 
As it can be seen from the graphs we are hitting the right distributions with the Metropolis step.

\vspace{0.1in}
\noindent
{\bf Example 2.} Here the length biased distribution is the mixture
$$
g(x)=0.25\,{\rm Ga}(x|\,2,1)+0.75\,{\rm Ga}(x|\,10,1) 
$$
We simulate a sample $y_{\rm mg}=(y_1,\ldots,y_n)$ 
of size $n=70$.
Similar results, as in the first example, are shown in Figure $2$, (a)--(c). 
For both estimators $\widetilde{g}_h(y;y_{\rm mg})$ and ${\widehat f}_{J,h}(y;y_{\rm mg})$ the bandwidth parameter 
has been calculated to $h=1.48$. For the average number of clusters, we estimate ${\cal C}_{\rm mg}=5.55$.
It is noted that the  Metropolis sampler  produces samples that are very close to the debiased  mixture
$f(x)=0.75\,{\rm Ga}(x|\,1,1)+0.25\,{\rm Ga}(x|\,9,1)$ depicted with a solid line in  $2$(c). 
 
\begin{figure*}
\includegraphics[width=160mm]{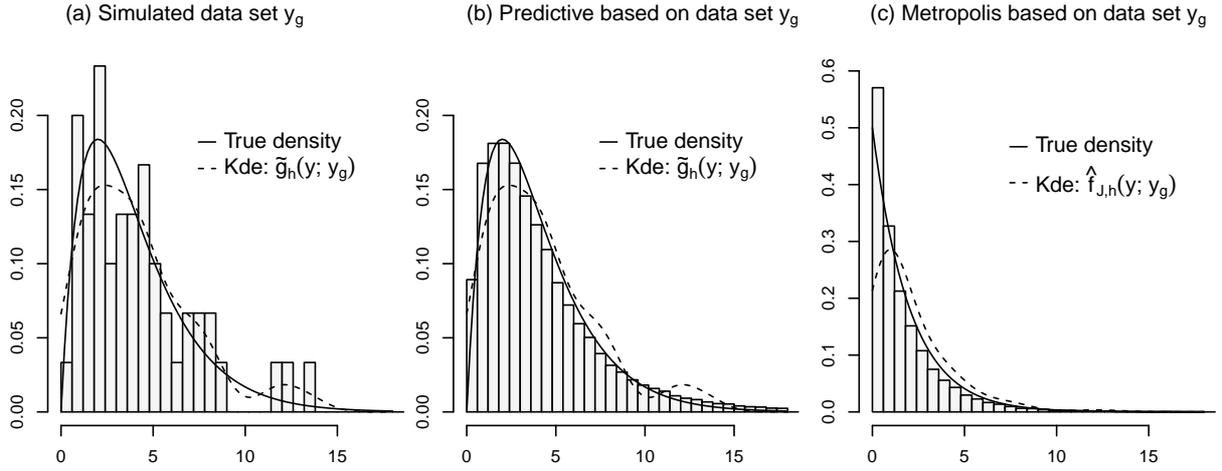}
\caption{ 
         Data set from  Ga$(2,0.5)$ of size $n=50$. 
         In all subfigures  the true densities are  depicted with a solid line and the kernel density estimates $\tilde{g}_h$ and $\hat{f}_{J,h}$  
         with a dashed line.}
\end{figure*}

\begin{figure*}
\includegraphics[width=160mm]{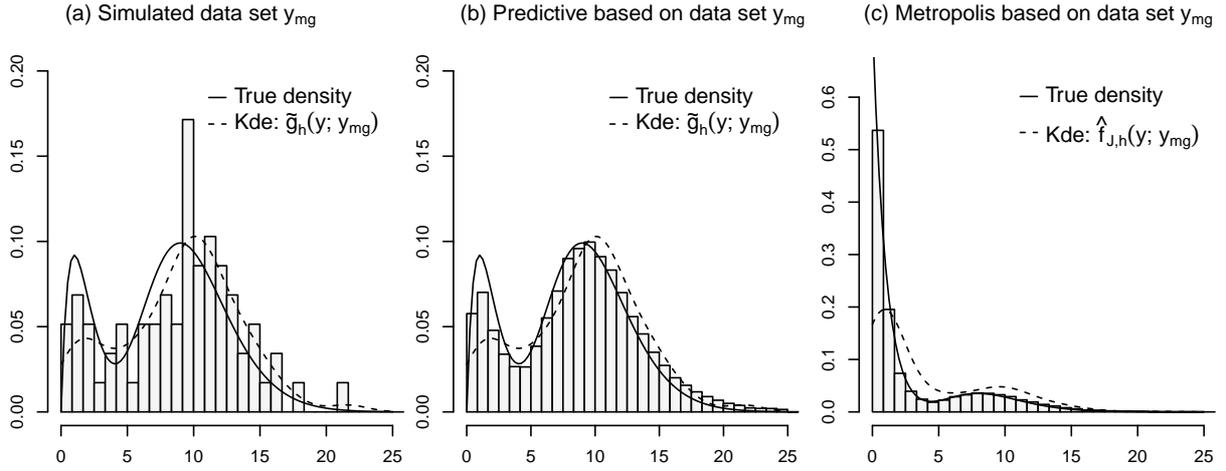}
\caption{
         Data set  from the mixture  $0.25\,{\rm Ga}(2,1)+0.75\,{\rm Ga}(10,1)$, $n=70$  .
         In all subfigures  the true densities are depicted with a solid line and  the kernel density estimates  $\tilde{g}_h$ and $\hat{f}_J$  
         with a dashed line.}
\end{figure*}
              
\begin{figure*}
\includegraphics[width=160mm]{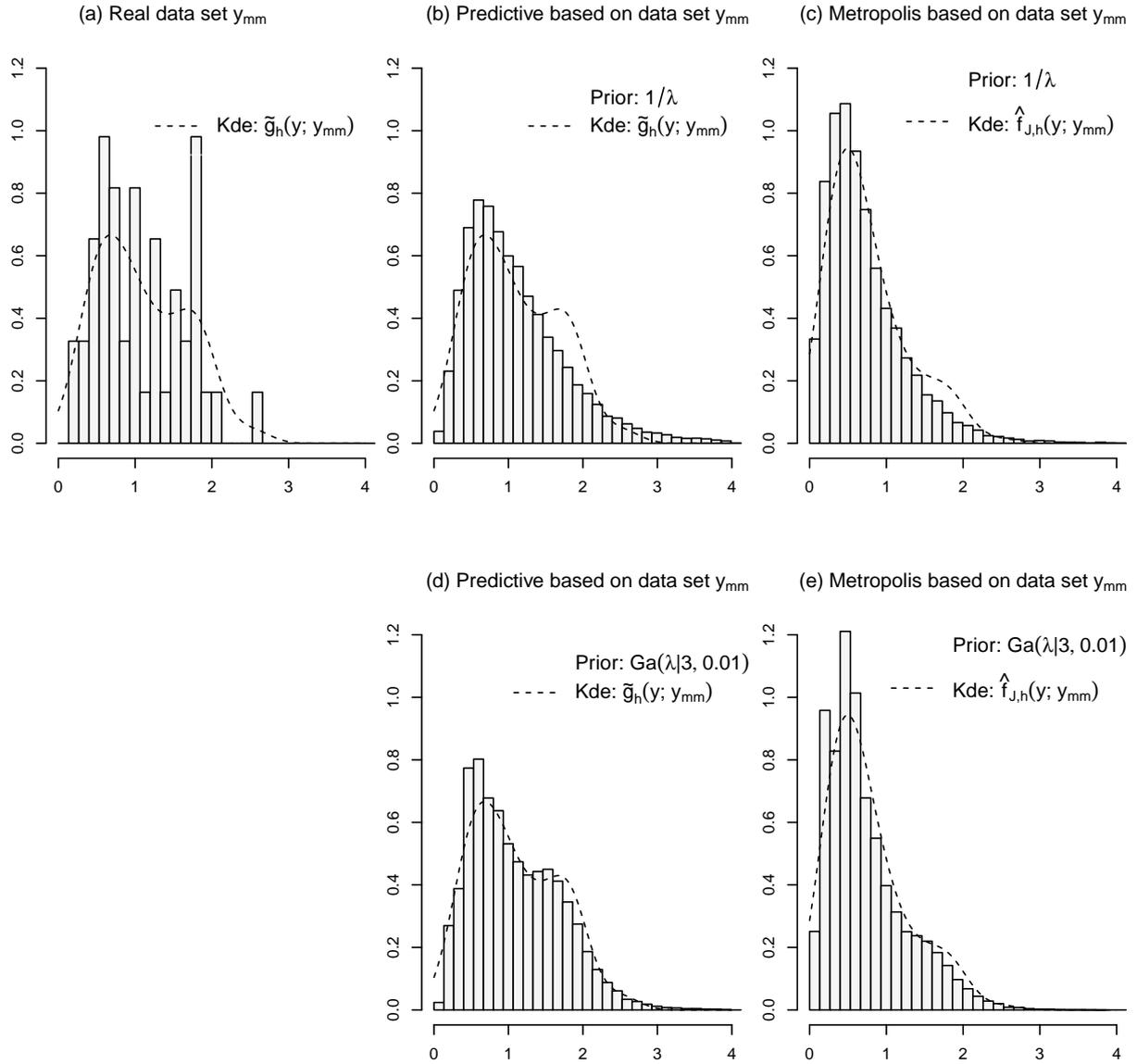}
\caption{
        Data set of size $n=46$  measuring the widths of shrubs.  Kernel density estimates 
        $\tilde{g}_h$ (classical) and $\hat{f}_{J,h}$ (indirect data) are depicted with dashed lines. 
        Top figures indicate  a noninformative prior specificaion while  bottom figures an informative one.
        Such  choice reproduces classical estimation results.}
\label{fig:3}
\end{figure*} 
\begin{figure*}
\includegraphics[width=160mm]{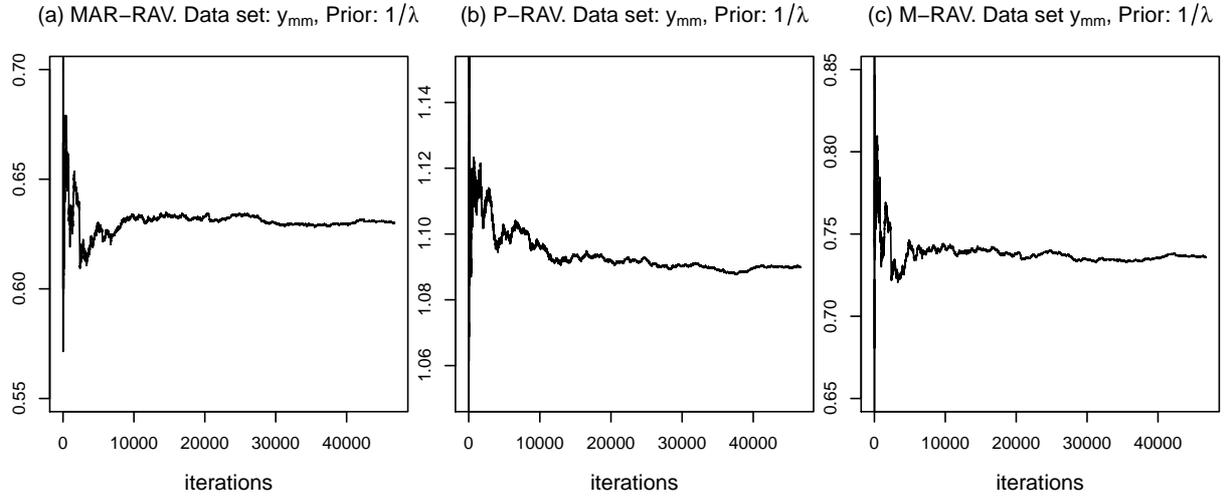}
\caption{ 
        (a) the Metropolis acceptance rate running average(MAR--RAV),(b) the predictive running average (P--RAV) and (c) the Metropolis sample running average (M--RAV).}
\label{figure2}
\end{figure*}    

\begin{figure*}
\includegraphics[width=160mm]{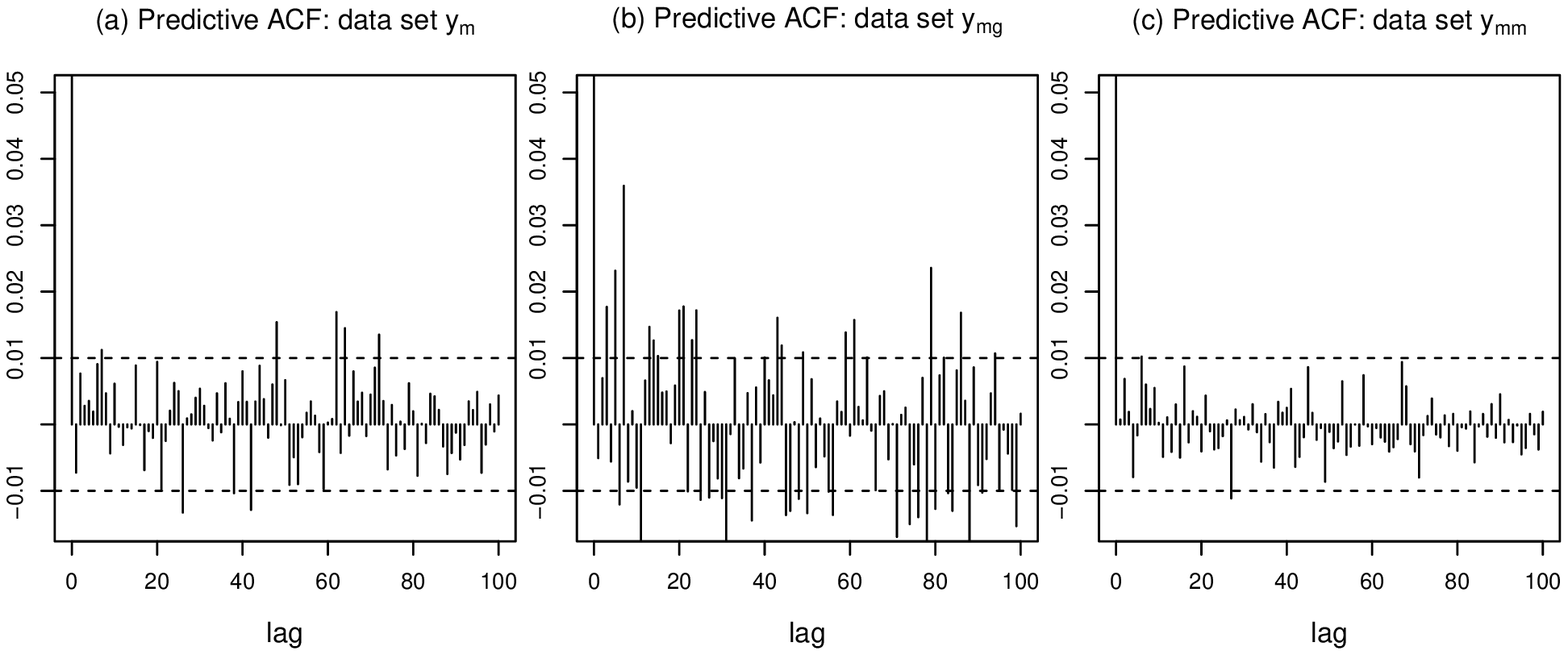}
\caption{
        The autocorrelation functions based on the posterior predictive observations. In (a),(b) the gamma and mixture of gammas synthetic data. In c)   the width of shrubs real data.}
\end{figure*}

\subsection{Real Data Example}

The data can be found in  Muttlak and McDonald (1990) and consist of $y_{\rm mm}=(y_1,\ldots,y_n)$,
$n = 46$, measurements  representing  the widths of shrubs obtained by line--transect sampling. 
In  this sampling method the probability of inclusion in the sample is proportional to the width 
of the  shrub  making it a case of length biased sampling.
A noninformative estimation is shown in Figure $3$ (a)-(c) with the same specifications as in (\ref{noninformative})  while in $3$(d), 3(e) we perform  a highly informative estimation with  
$\pi(\lambda)={\rm Ga}(\lambda|\,3,0.01)$.

The following results are presented  in Figures $3$ and $4$:

\begin{itemize}

\item $3$(a), $3$(b):  histograms of the length biased data set $y_{\rm mm}$ and  of  a sample from the posterior predictive $g(y_{n+1}|\,y_{\rm mm})$, respectively . In both subfigures the associated classical estimator $\widetilde{g}_h(y;y_{\rm mm})$ is depicted with a dashed line, for $h=0.23$. 
\item $3$(c): a histogram of the debiased data associated with the Metropolis chain estimator. Jones' estimator 
${\widehat f}_{J,h}(y;y_{\rm mm})$ is shown in dashed line, for the same bandwidth value. 

\item $3$(d), $3$(e): histograms of the posterior predictive and the Metropolis
sample, respectively, under the highly informative prior $\pi(\lambda)={\rm Ga}(\lambda|\,3,0.01)$, 
with superimposed classical density estimators.

\item $4$(a): the running acceptance rate of the Metropolis with jump distribution the posterior
predictive values from $g(y_{n+1}|\,y_{\rm mm})$ with an estimated value of about $0.62$.

\item $4$(b), $4$(c): running averages of the predictive and Metropolis samples respectively. 
\end{itemize}

Finally, in Figure 5 we provide the  autocorrelation function  as a function of lag, among the values of the  posterior predictive sample for the synthetic and real data sets, after a reasonable burn-in period.

\subsection{Remarks}

\begin{itemize}
\item Estimation  for the simulated data is nearly perfect and we get the best results 
for $ \pi(\lambda)\propto 1/\lambda $. As it is evident from subfigure $1$(c),  for the $\exp(0.5)$,  
the estimator ${\hat f}_{J,h}$ does not properly capture  the distributional features near the origin.
The same holds  true for the 'debiased' mixture density $0.75\,{\rm Ga}(x|\,1,1)+0.25\,{\rm Ga}(x|\,9,1)$, subfigure $2$(c).

\item For the real data set the $1/\lambda$  prior gives again the best results. Such a prior gives the largest average number of 
clusters among all noninformative specifications that were  examined.  The debiased $f$ density is close to ${\widehat f}_{J,h}$ 
though not exactly  the same. The difference comes from  a small area where the biased data have the group of observations 
$(1.85, 1.85, 1.86)$ that causes ${\widehat f}_{J,h}$ to produce an intense second mode. Excluding these $3$ data points Jones' estimator 
${\widehat f}_{J,h}$ becomes identical with ours. 

\item The highly informative specification 
$\lambda\sim{\rm Ga}(3,0.01)$ increases the average number of clusters from $4.03$ (noninformative estimation) to about $5.63$, thus the appearance of a second mode between $1.2$ and $2.0$, in $3$(d). From our numerical experiments it seems  that ${\widehat f}_{J,h}$  is 
"data hunting" in the sense that it  overestimates data sets and produces spurious modes. Our method performs 
better as it does not tend to overestimate, and at the same time has better properties near the origin.

\item When informative prior specifications are used they increase the average number of realized clusters and  the nonparametric estimates
tend to look more like Jones' type  estimates.
For example choices of $\lambda$ priors like Ga$(\alpha,0.01)$ with $\alpha\ge 2.5$ increase considerably  the average number of clusters 
and our real data estimates  in subfigures $3$(d) and 3(e) become nearly identical to ${\widehat f}_{J,h}$ .
\end{itemize}

\section{ Discussion}

In this paper we have described a novel approach to the Bayesian nonparametric modeling of a length bias sampling model. We directly tackle the length bias sampling distribution, from where the data arise, and this technique avoids the impossible situation of the normalizing constant if one decides to model the density of interest directly. This is legitimate modeling since only mild assumptions are made  on both densities, so we are free to  model  $g$ directly and choose  an appropriate  kernel with  the only condition that $\int_0^\infty x^{-1}k(x,\theta) dx<\infty$.

In a parametric set-up since $f$  is known up to a parameter $\theta$  modeling $g$ directly is not recommended, since to avoid a normalizing constant problem a model for $g$ would not result from the correct family for $f$.

We have also as part of the solution presented a Metropolis step to ``turn" the samples from $g(\cdot)$ into samples from $f(\cdot)$. A rejection sampler here would not work as the $1/y$ is unbounded.

The method we have proposed here should also be applicable to an arbitrary weight function $w(y)$, whereby samples are obtained from $g(y)$ and yet interest focuses on the density function $f(y)$, where the connection is provided by
$$g(y)=\frac{w(y)\,f(y)}{\int_0^\infty w(y)\,f(y)\,d y}.$$

Our estimator, besides being the first Bayesian kernel  density estimator for length biased   data, it was demonstrated that it performs at least as well and in some cases even better than its frequentist  counterpart.

\section *{Appendix: Asymptotics} 
In this section we assume that the posterior predictive sequence $(g_n)_{n\geq 1}$ is consistent in the sense that $d_1(g_n,g_0)\rightarrow 0$ a.s. as $n\rightarrow\infty$, where $g_0$ is the true density function generating the data and $d_1$ denotes the $L_1$ distance.
This would be a standard result in Bayesian nonparametric consistency involving mixture of Dirichlet  process models: see, for example, Lijoi et al. (2005), where sufficient conditions for the $L_1$ consistency are given.

The following theorem establishes a similar consistency result for the  debiased  density.

\vspace{0.1in}
\noindent
{\bf Theorem.} Let $f_n(y)\propto y^{-1}g_n(y)$ and $f_0(y)\propto y^{-1}g_0(y) $
denote  the  sequence of posterior predictive estimates  for the debiased density   and the true debiased density, respectively.
Then,  $d_1(f_n,f_0)\rightarrow 0$ a.s.

\vspace{0.1in}
\noindent
{\bf Proof.}  Let $$g_n(y)=\int \mbox{LN}(y|\mu,\sigma^2)\,dP_n(\mu,\sigma^2)$$
where $P_n$ is the posterior expectation of $P$, and for some $P_0$ it is that
$$g_0(y)=\int \mbox{LN}(y|\mu,\sigma^2)\,d P_0(\mu,\sigma^2).$$
The assumption of consistency also implies that $P_n$ converges weakly to $P_0$ with probability one.
This means for any continuous and bounded function $h(\mu,\sigma^2)$ of $(\mu,\sigma^2)$ we have the a.s. weak consistency of $P_n$ implies
$$\int h(\mu,\sigma^2)\,d P_n(\mu,\sigma^2)\rightarrow \int h(\mu,\sigma^2)\,d P_0(\mu,\sigma^2)\quad\mbox{a.s.}$$
and note that
$$\int_0^\infty y^{-1}\,\mbox{LN}(y|\mu,\sigma^2)\,d y=\exp\{-\mu+\sigma^2/2\}.$$

We now aim to show that these results imply the a.s. $L_1$ convergence of $f_n(y)$ to $f_0(y)$.
To this end, if we construct the prior so that for some constants $M$ and $\sigma^*$ it is that $\underline{\sigma}<\sigma<\overline{\sigma}$ and $|\mu|<M$, assuming $P_0$ puts all the mass on $[-M,+M]\times (\underline{\sigma}^2,\overline{\sigma}^{2})$, then from the definition of weak convergence we have that, with probability one,
\begin{equation*}
\begin{split} c_n=\int y^{-1} g_n(y)\,d y=\int \exp\{-\mu+\sigma^2/2\}\,d P_n(\mu,\sigma^2)\\\rightarrow c_0=
\int \exp\{-\mu+\sigma^2/2\}\,d P_0(\mu,\sigma^2).
\end{split} 
\end{equation*}

Also, with the conditions on $(\mu,\sigma^2)$,  we have
$$h_y(\mu,\sigma^2)=y^{-1}\mbox{LN}(y|\mu,\sigma^2)$$
is a bounded and continuous function of $(\mu,\sigma^2)$ for all $y>0$. 
Hence
\begin{equation*}
\begin{split}
\int y^{-1}\mbox{LN}(y|\mu,\sigma^2)\,d P_n(\mu,\sigma^2)\rightarrow \\\rightarrow
\int y^{-1}\mbox{LN}(y|\mu,\sigma^2)\,d P_0(\mu,\sigma^2)\quad\mbox{a.s.}\end{split}
\end{equation*}

pointwise for all $y>0$. Consequently, by Scheff\'e's theorem, we have
$$\int y^{-1}|g_n(y)-g_0(y)|\,d y\rightarrow 0\quad\mbox{a.s.}\,.$$
Now
$$\begin{array}{ll}
\int |f_n(y)-f_0(y)|\,d y \leq \int y^{-1}g_n(y)\,d y|c_n^{-1}-c_0^{-1}|+ &\\ +c_0^{-1}\int|y^{-1}g_n(y)-y^{-1}g_0(y)|\,d y\leq &  \\ 
\leq |1-c_n/c_0|+c_0^{-1}\int y^{-1}|g_n(y)-g_0(y)|\,d y &\\
\end{array} $$
and so
$$\int |f_n(y)-f_0(y)|\,d y\rightarrow 0\quad\mbox{a.s.}\,,$$
as required.



%
%

\section*{References}

\begin{description}

\item {\sc Adams, R. P.,  Murray, I. and MacKay, D.J.C.} 
The Gaussian process density sampler. {\sl Advances in Neural Information Processing Systems (NIPS)} {\bf 21}(2009).

\item {\sc Asgharian, M., M'Lan, C.E. and Wolfson, D.B.} 
Length--biased sampling with right--censoring: an unconditional approach.
{\sl Journal of the American Statistical Association} {\bf 97}, 201-209(2002).


\item   {\sc Cox, D.R.} 
Some sampling problems in technology.
In {\sl New Developments in Survey Sampling}. U.L. Johnson \& H. Smith eds. Wiley, New York(1969).



\item  {\sc Efron, B.}
Poisson overdispersion estimates based on the method of asymmetric maximum likelihood.
{\sl Journal of the American Statistical Association} {\bf 87}, 98--107(1986).

\item {\sc Escobar M. D. and West M.} Bayesian density estimation and inference using mixtures.
{\sl Journal of the American Statistical Association} {\bf 90}, 577--588(1995).


 \item{\sc Ferguson, T.S. }
A Bayesian analysis of some nonparametric problems.
{\sl Annals of Statistics} {\bf 1}, 209--230(1973).


\item  {\sc Jones, M.C.}  Kernel density estimation for length biased data.
{\sl Biometrika} {\bf 78}, 511--519(1991).

\item {\sc Kalli, M., Griffin, J. E. and Walker, S. G.}
Slice sampling mixture models. {\sl Statistics and Computing}  {\bf 21}, 93--105)(2011).

\item{\sc Lijoi, A., Pruenster, I. and Walker, S.G.} (2005). On consistency of nonparametric
normal mixtures for Bayesian density estimation. {\sl Journal
of the American Statistical Association} {\bf 100} , 1292--1296(2005).

\item {\sc Lo, A.Y.}
On a class of Bayesian nonparametric estimates I. Density estimates.
{\sl Annals of Statistics} {\bf 12}, 351--357(1984).

\item  {\sc Murray, I., Ghahramani, Z. and MacKay, D.J.C.}
MCMC for doubly intractable
distributions. {\sl Proceedings of the 22nd Annual Conference on Uncertainty in Artificial
Intelligence (UAI)}, 359--366(2006).

\item{\sc Muttlak, H.A. and McDonald, L.L.} 
Ranked set sampling with size-biased probability of selection.
{\sl Biometrics} {\bf 46}, 435--445(1990).

\item {\sc Patil, G.P. } Weighted distributions.
{\sl Encyclopedia of Environmetrics} (ISBN 0471 899976)) {\bf 4}, 2369--2377(2002).

\item  {\sc Patil, G.P. and Rao, C.R.} 
Weighted distributions and size-biased sampling with applications to
wildlife populations and human families.
{\sl Biometrics} {\bf 34}, 179--189(1978).



\item{\sc Sethuraman, J. and Tiwari, T.C.}
Convergence of Dirichlet measures and the interpretation of their parameter.
{\sl Statistical Decision Theory and Related Topics III} \textbf{2}, 305--315(1982).

 \item{\sc Sethuraman, J.} A constructive definition of Dirichlet priors.
{\sl Statistica Sinica} {\bf 4}, 639--650(1994).

\item{\sc Sheather, S.J.  and Jones, M.C.} A Reliable Data-based Bandwidth Selection Method
for Kernel Density Estimation. {\sl Journal
of the Royal Statistical Society  B} {\bf 12}, 683--690(1991).

\item{\sc Tokdar, S.T.} Towards a faster implementation of density estimation with logistic  Gaussian process prior. {\sl Journal
of Computational and Graphical Statistics} {\bf 16} , 633--655(2007).

\item{\sc Vardi, Y.}  Nonparametric estimation in the presence of length bias.
{\sl The Annals of Statistics} {\bf 10}, 616--620(1982).

\item {\sc Walker, S.G.}  Sampling the Dirichlet mixture model with slices.
{\sl Communications in Statistics} {\bf 36}, 45--54(2007).

\item{\sc Walker, S.G.} Posterior sampling when the normalizing constant is unknown.
{\sl Communications in Statistics} {\bf 40}, 784--792(2011).

\end{description}

\end{document}